\documentclass[a4paper]{amsart}
\usepackage{graphics}
\usepackage{amsmath}
\usepackage{latexsym}
\usepackage{amssymb}
\usepackage[all]{xy}
\xyoption{matrix}
\xyoption{arrow}

\begin{document}

\newtheorem{lem}{Lemma}[section]
\newtheorem{prop}[lem]{Proposition}
\newtheorem{cor}[lem]{Corollary}
\newtheorem{thm}[lem]{Theorem}
\newtheorem*{thmA}{Theorem}
\newtheorem*{thmB}{Theorem B}
\newtheorem{rem}[lem]{Remark}
\newtheorem{defin}[lem]{Definition}

\newcommand{\Z}{\operatorname{\mathbb Z}\nolimits}
\newcommand{\Q}{\operatorname{\mathbb Q}\nolimits}
\newcommand{\Ext}{\operatorname{Ext}\nolimits}
\newcommand{\Hom}{\operatorname{Hom}\nolimits}
\renewcommand{\mod}{\operatorname{mod}\nolimits}
\renewcommand{\dim}{\operatorname{dim}\nolimits}
\newcommand{\End}{\operatorname{End}\nolimits}
\newcommand{\op}{\operatorname{op}\nolimits}
\newcommand{\C}{\operatorname{\mathcal C}\nolimits}
\newcommand{\M}{\operatorname{\mathcal M}\nolimits}
\renewcommand{\P}{\operatorname{\mathcal P}\nolimits}

\title[Cluster algebras of extended Dynkin type]
{Cluster algebras associated with extended Dynkin quivers} 

\author[Buan]{Aslak Bakke Buan}
\address{Institutt for matematiske fag\\
Norges teknisk-naturvitenskapelige universitet\\
N-7491 Trondheim\\
Norway}
\email{aslakb@math.ntnu.no}

\author[Reiten]{Idun Reiten}
\address{Institutt for matematiske fag\\
Norges teknisk-naturvitenskapelige universitet\\
N-7491 Trondheim\\
Norway}
\email{idunr@math.ntnu.no}

\maketitle

\section{Introduction}

Fomin and Zelevinsky \cite{fz1} have defined cluster algebras, 
and developed an interesting and influential theory about this class of algebras.
We deal here with a special case of cluster algebras.
Let $B$ be an $n \times n$ integral skew symmetric matrix,
or equivalently a  finite quiver $Q_B$ with $n$
vertices and no loops or oriented cycles of length two. Let
$\underline{u} = \{u_1, \dots, u_n \}$ be a
transcendence basis for $F= \Q(x_1, \dots, x_n)$.
Then there is a cluster algebra associated to the pair $(B, \underline{u})$.

An essential ingredient in the definition of cluster algebras is the concept of mutation of matrices,
or equivalently mutation of quivers. For each row $i=1,\dots,n$ in
the skew-symmetric matrix $B$, there is an associated skew-symmetric matrix $\mu_i(B)$.
If the quiver $Q_B$ is a Dynkin quiver, it follows from \cite{fz1} that there is only a finite number of
non-isomorphic quivers obtained from $Q_B$ by sequences of mutations. The purpose of this paper is to
show that when $Q_B$ is a connected quiver with no oriented cycles, then there is only a finite number
of quivers obtained by sequences of mutations if and only if
$Q$ is Dynkin or extended Dynkin, or has at most two vertices. 

The first link from cluster algebras to tilting theory for finite dimensional algebras was
discovered by Marsh, Reineke and Zelevinsky \cite{mrz}. This inspired the invention
of cluster categories and cluster-tilted algebras \cite{bmrrt, bmr1, bmr2}.
Our main result is obtained
as an application of this work. 
It answers a question by A. Seven, who has showed one 
implication (that finite mutation type implies Dynkin, extended Dynkin or
rank two) using different methods \cite{s1}.

We would like to thank Otto Kerner for very helpful conversations.

\section{Background}

In this section we collect some background material on finite dimensional hereditary
algebras, cluster algebras and on cluster categories and cluster-tilted algebras.

Let $K$ be an algebraically closed field, and $Q$ a finite connected quiver with no oriented cycles.
The associated path algebra $KQ$ is a finite dimensional hereditary algebra. For such algebras there is a well developed
module theory, see \cite{r,ars}. 
In particular there is a nice characterization of finite, tame and wild representation type
in this case.

\begin{thm}
Let $KQ$ be a finite-dimensional hereditary $K$-algebra. Then we have the following.
\begin{itemize}
\item[(a)]{$KQ$ is of finite representation type if and only if the underlying
graph of $Q$ is Dynkin.}
\item[(b)]{$KQ$ is of tame representation type if and only if the underlying graph is extended Dynkin.}
\end{itemize}
\end{thm}

In the case of infinite type there are three kinds of indecomposable modules: preprojective, preinjective and regular.
The preprojective ones are those of the form $\tau^{-i}P$ for $i \geq 0$, where $P$ 
is projective and $\tau$ is the AR-translation.
The preinjective ones are those of the form $\tau^i I$ for $i \geq 0$, where $I$ 
is injective, and the regular ones are the rest. 
All the preprojective (preinjective) indecomposables are in the same component of the AR-quiver, this is called
the preprojective (preinjective) component. 
For the regular ones, the components of the AR-quiver are tubes in the tame case and of 
the form $\Z A_{\infty}$ in the wild case.
In each case one has the concept of {\em quasi-length} for regular modules. Informally,
the quasi-length of an indecomposable regular module is a measure of how far the module
is from the border of the component. See \cite{ker} for details.

A $KQ$-module $X$ is said to be exceptional if $\Ext^1_{KQ}(X,X) = 0$. There is the following
useful information about such modules, where part (b) is from \cite{h}.
  
\begin{prop}\label{quasi}
Let $KQ$ be a finite dimensional hereditary algebra, and $X$ an indecomposable $KQ$-module.
\begin{itemize}
\item[(a)]{If $X$ is preprojective or preinjective, then
$X$ is exceptional.}
\item[(b)]{If $X$ is regular and exceptional, then $X$ has quasi-length at most $n-2$,
where $n$ denotes the number of vertices in the quiver.}
\end{itemize}
\end{prop}

A $KQ$-module $T$ is a {\em tilting module}
if $T$ is exceptional and there is an exact sequence $0 \to KQ \to T_0 \to T_1 \to 0$,
with $T_0$ and $T_1$ direct summands of finite direct sums of copies of $T$. We have the following
\cite{hr, b}.

\begin{prop}
\item[(a)]{If $T$ is a tilting $KQ$-module, then the number of non-isomorphic indecomposable direct summands
of $T$ is equal to the number of vertices of $Q$.}
\item[(b)]{Any exceptional $KQ$-module can be extended to a tilting module.
}
\end{prop}

We next give some informal background on cluster algebras \cite{fz1}. Let $n$ 
be a positive integer and $B = (b_{ij})$ a skew symmetric integral $n \times n$-matrix.
Associated with $B$ is a quiver $Q = Q_B$ with vertices $1, \dots ,n$ corresponding to the
rows of the matrix. 

Then $Q$ has no loops and no oriented cycles of length two.
If $b_{ij} > 0$, there are $b_{ij}$ arrows from $i$ to $j$.
Let $\mathbb{F} = \Q(x_1, \dots , x_n)$, where the $x_i$ are indeterminates.
Start with a so-called seed $(\underline{x},B)$ where $\underline{x} = \{x_1, x_2, \dots x_n \}$ 
is a transcendence basis and $B$ an integral skew-symmetric $n \times n$-matrix.
For each index  $k = 1, \dots,n$ we define a new seed $\mu_k(\underline{x},B) = (\underline{x'},B')$
as follows.
First define a new element $x_k' \in\mathbb{F}$ via the relation:
$x_k' x_k =\prod_{x\in\underline{x},b_{ik}>0}x^{b_{ik}}+
\prod_{x\in\underline{x},b_{ik}<0}x^{-b_{ik}}.$
Here, we say that $x_k,x_k'$ form an {\em exchange pair}.
Let $\underline{x}'=\underline{x}\cup \{x_k'\}\setminus \{x_k\}$, this is a new
transcendence basis of $\mathbb{F}$. Let $B' = (b'_{ij})$ be the
{\em mutation} of the matrix $B$ in direction $k$ (as defined in~\cite{fz1}), that is
$$b'_{ij}=\left\{
\begin{array}{ll}
-b_{ij} & \mbox{if\ } i=k \mbox{\ or\ }j=k, \\
b_{ij}+\frac{1}{2}(|b_{ik}|b_{kj}+b_{ik}|b_{kj}|) & \mbox{otherwise.}
\end{array}
\right. $$
The pair $(\underline{x'},B')$ is called the {\em mutation} of the seed
$(\underline{x},B)$ in direction
$k$. 
The transcendence bases obtained this way are called {\em clusters}, and the
elements in the clusters are called {\em cluster variables.}
The associated cluster algebra is the subalgebra of $F$ generated by the cluster 
variables. The case where one of the matrices occurring in a seed corresponds to
a quiver without oriented cycles gives by definition an {\em acyclic cluster algebra} \cite{bfz}.

Note that the mutation operation is defined on skew-symmetric integral matrices. Such matrices correspond to 
finite quivers with no loops or oriented cycles of length two. We let the 
{\em mutation class} of such a quiver $Q$ denote the set of all quivers which can be reached
from $Q$ by a finite number of mutations.

A major result is the following \cite{fz2}, where a cluster algebra is said to be of finite type if
there are only a finite number of seeds.

\begin{thm} 
The following are equivalent for a cluster algebra.
\begin{itemize}
\item[(a)]{There is only a finite number of seeds.}
\item[(b)]{There is only a finite number of clusters.}
\item[(c)]{One of the matrices occurring in a seed is associated with a Dynkin quiver.}
\end{itemize}
\end{thm}

In \cite{s2}, Seven gives a list of cluster algebras of socalled minimal infinite type. This can 
be used to determine if a given cluster algebra is of finite type.
If a cluster algebra is of finite type, there is also a finite number of associated non-isomorphic quivers.
But the converse is not true, and the aim of this paper is to characterize the class of
acyclic cluster algebras which have only a finite number of associated non-isomorphic quivers.

The cluster categories associated with finite dimensional hereditary algebras were
introduced and investigated in \cite{bmrrt},
in trying to model some of the concepts from cluster algebras in a categorical/module theoretical way.
Let as before $K$ be an algebraically closed field,
and $Q$ a finite quiver without oriented cycles, and $H= KQ$ the associated path algebra. Then the associated 
cluster category is by definition $\C_H = D^b(H)/ \tau^{-1}[1]$, where $D^b(H)$ is the bounded derived
category of the finitely generated $H$-modules and $[1]$ the shift functor in $D^b(H)$.
Then $\C$ is a triangulated category \cite{k}. (Cluster-)tilting objects in $\C_H$ are defined to be basic 
objects $T$ such that $\Ext^1_{\C_H}(T,T) = 0$, and such that $T$ is maximal with respect to 
this property. 
The tilting objects in $\C_H$ are the analogs of clusters, and the indecomposable objects
in $\C_H$ are the analogs of cluster variables.
The tilting objects in $\C_H$ are the objects induced (via a canonical functor $\mod H' \to \C_H$) by a tilting
module over some hereditary algebra $H'$ derived equivalent to $H$.
Their endomorphism algebras are the cluster-tilted algebras investigated in \cite{bmr1, bmr2}.

In the next section we use the theory of cluster categories and cluster-tilted algebras, along with 
results on hereditary algebras, to prove our main result. A crucial role in this investigation
is played by the following \cite{bmr2}.

\begin{thm}\label{thesame}
Let $Q$ be a finite quiver without oriented cycles. Then the mutation class of $Q$ 
coincides with the set of quivers of cluster-tilted algebras determined by $KQ$.
\end{thm}

\section{Main result}

This section is devoted to proving when 
a quiver without oriented cycles has a finite mutation class, and thus proving when acyclic cluster algebras
are defined using only a finite number of non-isomorphic quivers. This is obtained as 
a direct consequence of the following result on cluster-tilted algebras.

\begin{thm}\label{counting}
Let $K$ be an algebraically closed field and $Q$ a connected finite quiver without
oriented cycles. The the following are equivalent for $H=KQ$.
\begin{itemize}
\item[(a)]{There is only a finite number of (basic) cluster-tilted algebras associated with $H$,
up to isomorphism.}
\item[(b)]{There is only a finite number of quivers, up to isomorphism,
occurring for cluster-tilted algebras associated with $H$.}
\item[(c)]{$H$ is of finite or tame representation type, or has at most two non-isomorphic
simple modules.}
\end{itemize}
\end{thm}

\begin{proof}
(a) implies (b): This is obvious.

\noindent (c) implies (a): If $H$ is of finite representation type,
the statement in (a) is obvious, since there are only a finite number of 
non-isomorphic tilting objects in this case.

For the remaining cases, the following easily seen fact is useful.

\begin{lem}
For any integer $n$ , and any object $X$ in $\C_H$, we have $\End_{\C_H}(X) \simeq \End_{\C_H}( \tau^n X)$.
\end{lem}

\begin{proof}
This follows directly from the fact that $\tau \colon \C_H \to \C_H$ is an autoequivalence.
\end{proof}

Assume next that there are at most two (non-isomorphic) simple $H$-modules. If there is only one simple
module we have $KQ \simeq K$. So we can assume there are two simples. Then there
is no indecomposable regular module which is exceptional, by Proposition \ref{quasi}.
Up to $\tau$-shift in the cluster category $\C_H$ we can assume that a tilting object in $\C_H$
is induced by a tilting $H$-module $T$ with an indecomposable projective direct summand
$P$. Then the almost complete tilting object $P$ has exactly two complements in $\C_H$ by \cite{bmrrt}. 
Since there are two possible choices for 
$P$, we obtain that there are at most four quivers of cluster-tilted algebras.

Assume now that $H$ is tame.
We show that there is a finite set $\M = \{M_i \}_{i \in I}$, such that for any tilting object $T$ in $\C_H$,
there is an integer $j$ and an $i \in I$, such that $T \simeq  \tau^j M_i$.
We need the following.

\begin{lem}\label{finite}
For each indecomposable projective $H$-module $P$,  
there is a finite set of objects $X$ in $\C_H$ with $\Ext^1_{\C}(X,P)= 0$. 
\end{lem}

\begin{proof}
All indecomposable objects in $\C_H$ are either isomorphic to $H$-modules (via the canonical map $\mod H \to \C_H$)
or isomorphic to objects of the form $P'[1]$, where $P'$ is indecomposable projective. 
It is therefore clearly sufficient to show that there is only a finite number of $H$-modules $X$
with $\Ext^1_{\C}(X,P)= 0$.
We have $\Ext^1_H(X,P) \simeq D\Hom_H(P, \tau X)$, where $D$ denotes the duality for $\mod H$. 
Let $S$ be the simple top of $P$. 
Then we have that $\Hom_H(P,Y) = 0$ if and only if $S$ is not a composition
factor of $Y$. For tame 
hereditary algebras, there is only a finite number of indecomposable modules which do not have $S$
as a composition factor. Thus, the claim follows.
\end{proof}

It follows from Lemma \ref{finite} that there is only a finite number of tilting objects with 
an indecomposable projective $H$-module $P$ as a direct summand,
and thus the set $\M$ of tilting objects in $\C_H$ with a projective $H$-module as a direct summand is finite.
It is known \cite{hr2} that $T$ has at least one direct summand which is preprojective or preinjective.
This finishes the proof that (c) implies (a). \\
\\
\noindent (b) implies (c):
Assume that there is only a finite number of quivers (up to isomorphism) occurring for cluster-tilted
algebras associated with $H$. Let $n$ be the number of vertices, and let $u$ be the maximum number
of arrows between two vertices in this finite set of quivers. From the following lemma we can then conclude that 
$\dim \Hom_{\C_H}(T_1,T_2) < m = u^{2n}$ for any two indecomposable direct summands $T_1, T_2$ of a tilting object $T$.
 
\begin{lem}
If a path in the quiver of a cluster-tilted algebra goes through two oriented cycles, then it is zero.
\end{lem}

\begin{proof}
Let $T  = \amalg_{i=1^n} T_i$ be a basic tilting module for a hereditary algebra $H$, with corresponding cluster
category $\C_H$. There are no oriented cycles in a tilted algebra. Thus, we know that one of the maps in an oriented
cycle for the cluster-tilted algebra is a map which lifts to a map of the
form $T_{i_1} \to F(T_{i_2}) = \tau^{-1}T_{i_2}[1]$ in $D^b(H)$.
Using that any map $T_a \to T_b[2]$ is zero in $D^b(H)$ and thus is zero also in $\C_H$, the claim is proved.
\end{proof}

Assume now that $H$ is wild with at least three non-isomorphic simple $H$-modules.
Then there is a vertex $e$ in the quiver of $H$ such that $A = H/HeH$ is of infinite representation type.
By possibly replacing $H$ by a derived equivalent hereditary algebra, we can assume that $e$ is a source
in the quiver. This can be seen as follows.
Consider the 
set $\P$ of paths $\rho$ in $Q$ (where $H= KQ$) with the following properties:
\begin{itemize}
\item[-]{The path $\rho$ ends in $e$.}
\item[-]{The path $\rho$ is not a proper subpath of a path $\rho'$ ending in $e$.}
\end{itemize}
Let $m_H$ be the sum of the lengths of the paths in $\P$. Assume $m_H > 0$. Then 
there is at least one path in $\P$, assume this paths starts in $f$. Then $f$ is a source in
the quiver $Q$. Consider the quiver $Q'$ obtained by reversing all arrows starting in
$f$ (this is a Bernstein-Gelfand-Ponomarev reflection \cite{bgp}).
Then $KQ' = H'$ is derived equivalent to $H$, and $m_{H'} < m_H$. 
Repeating this a finite number of times (at most $m_H$) we obtain
an algebra $H''$ derived equivalent to $H$, and with $m_{H''} = 0$. That is, $e$ is a source
for this algebra.

Choose $A'$ to be a connected summand of $A$ of infinite representation type. Let $S$ be the simple
$H$-module corresponding to $e$, and $P$ its projective cover. Choose a vertex $e'$ in the quiver 
of $A'$ which is a neighbour of $e$ in the quiver of $H$. Denote by $P'$ the corresponding indecomposable 
projective $A'$-module and let $S'$ be the simple top of $P'$. 
Let $X_i = \tau_{A'}^{-i}P'$ for $i \geq 0$. Then $\Ext_{A'}^1(X_i,X_i) = 0$ and hence
$\Ext_H^1(X_i,X_i) = 0$. Furthermore, $\Ext_H^1(P,\tau_H^{-1}X_i) = 0$, and  $\Ext_H^1(\tau_H^{-1}X_i,P) 
\simeq D\Hom_H(P,X_i)= 0$, since $S$ is not a composition factor of the $A'$-module $X_i$.
Hence $P \amalg \tau^{-1}_H X_i$ is an exceptional $H$-module for all $i \geq 0$. For each
$i \geq 0$ choose a tilting $H$-module $Y_i$ extending $P \amalg \tau^{-1}_H X_i$. 

There is some $i_0$
such that the simple $A'$-module $S'$ has multiplicity at least $m$ in $X_{i_0}$ by \cite{dr}.
Let $C$ denote the Coxeter
transformation, defined on $K_0(\mod H)$ (see \cite{bgp,ap}).
Then $C^{-1}[X_{i_0}] = [\tau^{-1}_H(X_{i_0})]$, where $[X]$ denotes the image of the $H$-module
$X$ in the Grothendieck group $K_0(\mod H)$. In particular, the multiplicity of $S'$ in $\tau^{-1}_H(X_{i_0})$
is the coefficient of $[S']$ when $C^{-1}[X_{i_0}]$ is expressed in terms of the basis coming from the 
simple $H$-modules. Since $e$ is a source in the quiver of $H$, we can compute this coefficient as follows,
by using the definition of the Coxeter transformation as a composition of reflections.
Let $e_1, \dots, e_t$ 
be the neighbours of $e$ in the quiver of $H$, with
$r_i$ arrows from $e$ to $e_i$, and let $m_i$ be the multiplicity of the corresponding simple module
$S_i$ in $X_{i_0}$. Then $\sum_{i=1}^t r_i m_i$ gives the multiplicity of
$S$ in $\tau^{-1}_H(X_{i_0})$, and this is at least $m$, since $S'$ is one of the $S_i$.

Hence $\dim \Hom_H(P,\tau^{-1}_H X_{i_0}) \geq m$, which gives a contradiction to the choice of $m$, since $P$
and $\tau^{-1}_H X_{i_0}$ are direct summands of the same tilting object.
\end{proof}

This has the following consequence for cluster algebras.

\begin{thm}
Let $A$ be an acyclic cluster algebra defined by a skew symmetric integral matrix, or equivalently
by a quiver $Q$. Then the following are equivalent.
\begin{itemize}
\item[(a)]{The mutation class of $Q$ is finite.}
\item[(b)]{$Q$ is either Dynkin, extended Dynkin or has at most two vertices.}
\end{itemize}
\end{thm}

\begin{proof}
This is direct consequence of Theorems \ref{thesame} and \ref{counting}.
\end{proof}

\bibliographystyle{plain}

\end{document}